\documentclass[10pt]{article}
\usepackage{amsmath,amssymb,amsfonts,latexsym,times,mathptmx}
\DeclareSymbolFont{symbolsYF}{OMS}{ptm}{m}{n}
\DeclareSymbolFontAlphabet{\mathcal}{symbolsYF}

\newtheorem{theorem}{Theorem}[section]
\newtheorem{proposition}{Proposition}

\newtheorem{lemma}[theorem]{Lemma}

\newtheorem{problem}[theorem]{Problem}

\title{Sparseness of $4$-cycle systems}

\author{Yuichiro Fujiwara
\thanks{Yuichiro Fujiwara is with the Graduate School of System and Information Engineering,
University of Tsukuba, Tsukuba, Ibaraki, Japan (e-mail: yuichiro@sk.tsukuba.ac.jp)}\and
Shung-Liang Wu
\thanks{Shung-Liang Wu is with National United University, Miaoli, Taiwan, R.O.C.}\and
Hung-Lin Fu
\thanks{Hung-Lin Fu is with Department of Applied Mathematics, National Chiao Tung University, Hsin Chu, Taiwan, R.O.C.}}

\date{}

\begin{document}
\maketitle

\begin{abstract}
\noindent
An avoidance problem of configurations in $4$-cycle systems is investigated
by generalizing the notion of sparseness,
which is originally from Erd\H{o}s'
$r$-sparse conjecture on Steiner triple systems.
A $4$-cycle system of order $v$, $4$CS$(v)$, is said to be $r$-sparse
if for every integer $j$ satisfying $2\leq j\leq r$
it contains no configurations consisting of $j$ $4$-cycles
whose union contains precisely $j+3$ vertices.
If an $r$-sparse $4$CS$(v)$ is also free from copies of a configuration on two $4$-cycles
sharing a diagonal, called the double-diamond, we say it is strictly $r$-sparse.
In this paper, we show that for every admissible order $v$ there exists a strictly $4$-sparse
$4$CS$(v)$. We also prove that for any positive integer $r \geq 2$ and sufficiently large
integer $v$ there exists a constant number $c$ such that there exists a strictly
$r$-sparse $4$-cycle packing of order $v$ with $c\cdot v^2$ $4$-cycles.
\end{abstract}

\begin{center}
\end{center}

\section{Introduction}
A $4$-{\it cycle system} of {\it order} $v$, denoted by $4$CS$(v)$, is an ordered pair $(V,{\mathcal C})$,
where $V=V(K_v)$, the vertex set of the complete graph $K_v$, and ${\mathcal C}$ is a collection
of edge-disjoint cycles of length four whose edges partition the edge set of the complete graph.
It is well-known that a necessary and sufficient condition for the existence of a $4$CS$(v)$
is that $v \equiv 1$ (mod $8$) (see, for example, Rodger \cite{R}).
Such orders are said to be {\it admissible}.
Following the usual terminology of cycle systems, we call a cycle of length four a $4$-{\it cycle}.

A $4$-cycle system is a natural generalization of the classical combinatorial design
called a Steiner triple system, briefly STS,
since an STS is just an edge-disjoint decomposition of a complete graph into triangles.
A Steiner triple system of order $v$ exists if and only if $v \equiv 1, 3$ (mod 6).
In other words, the set of all admissible orders of an STS consists of all the positive integers
$v \equiv 1, 3$ (mod 6).

As is the case with Steiner triple systems, various properties which may appear in
a $4$-cycle system
have also been studied (see, for example, Mishima and Fu \cite{MF} and references therein).
Such properties of cycle systems have also been investigated as a special graph design (see, for example, Jimbo and Kuriki \cite{JK}).
Among many characteristics of STSs, the numbers of occurrences of particular substructures
have been of interest to various areas (see Colbourn and Rosa \cite{CR}).
In the current paper, we consider an extreme case for $4$CSs,
namely, avoidance of particular configurations.
We first recall a long-standing conjecture on STSs posed by Erd\H{o}s.

A $(k,l)$-{\it configuration} in an STS is a set of $l$ triangles
whose union contains precisely $k$ vertices.
In 1973, Erd\H{o}s \cite{E} conjectured that for every integer $r \geq 4$, there exists
$v_0(r)$ such that if $v > v_0(r)$ and if $v$ is admissible, then there exists a Steiner
triple system of order $v$
with the property that it contains no $(j+2,j)$-configurations
for any $j$ satisfying $2\leq j\leq r$. Such an STS is said to be $r$-{\it sparse}.
Many results on the $r$-sparse conjecture and related problems have been since developed.
In particular,
after major progress due to Ling et al.\ \cite{LCGG} and earlier development found
in their references,
the simplest case when $r=4$, as it is sometimes called the anti-Pasch
conjecture, was eventually settled in the affirmative by Grannell et al.\ \cite{GGW}.

\begin{theorem}\label{anti-Pasch}{\rm\bf{(Grannell, Griggs and Whitehead \cite{GGW})}}
There exists a $4$-sparse Steiner triple system of order $v$
if and only if $v \equiv 1, 3$ {\rm (mod $6$)} and $v \not= 7, 13$.
\end{theorem}

As far as the authors are aware, the $r$-sparse conjecture for $r
\geq 5$ is still unsettled. In fact, no $7$-sparse STS is realized
for $v>3$. Very recent results on sparseness and related problems
are found in a series of papers: Forbes et al.\ \cite{FGG}, Wolf
\cite{W1,W2} and the first author \cite{F1,F2,F3,F4}. For general
background on configurations and sparseness in triple systems, the
interested reader is referred to Colbourn and Rosa \cite{CR}.

With regard to $4$-cycle systems, the relating result is due to
Bryant et al.\ \cite{BGGM}, who investigated the numbers of
occurrences of configurations consisting of two $4$-cycles. They
presented a formula for the number of occurrences of such
configurations and studied avoidance and maximizing problems.

Our primary focus in the current paper is on existence of $4$-cycle systems which are ``sparse"
in the sense that they do not contain configurations
that consist of many $4$-cycles on a small number of vertices in relative terms.
In this sense, for a given integer $w \leq v$ the ``densest" configurations on $w$ vertices
in a $4$CS$(v)$ are ones that contain as many $4$-cycles as possible.
In terms of combinatorial design theory,
such a configuration is said to be a
maximum $4$-cycle packing of order $w$.
More formally, a {\it $4$-cycle packing} of {\it order} $w$ is
an ordered pair $(W,{\mathcal D})$ such that $|W|=w$ and ${\mathcal D}$
is a set of $4$-cycles sharing no common edges, where vertices of a $4$-cycle
in ${\mathcal D}$ are elements of $W$.
A $4$-cycle packing is said to be ${\it maximum}$
if no other $4$-cycle packing of the same order contains a larger number of $4$-cycles.
Obviously, if $w$ is admissible, a maximum $4$-cycle packing of order $w$
is just a $4$CS$(w)$.

The term $(k,l)$-configuration will also be used for substructures in $4$CSs
and is defined as a set of $l$ $4$-cycles
on precisely $k$ vertices where no pair of distinct $4$-cycles share the same edge.
We denote the set of vertices in a configuration ${\mathcal A}$ by $V({\mathcal A})$.
Two configurations ${\mathcal A}$ and ${\mathcal B}$ are said to be {\it isomorphic},
denoted as ${\mathcal A} \cong {\mathcal B}$, if there exists a bijection
$\phi : V({\mathcal A}) \rightarrow V({\mathcal B})$ such that for each $4$-cycle
$C \in {\mathcal A}$, the image $\phi(C)$ is a $4$-cycle in ${\mathcal B}$.

In the case of STSs,
sparseness is measured by lack of $(j+2,j)$-configurations; one of reasons may be that
they are possibly avoidable and form the essential portions of dense configurations
(see Forbes, Grannell and Griggs \cite{FGG}).
Based on the following proposition and subsequent observation on $(j+3,j)$-configurations,
we propose an avoidance problem similar to the $r$-sparse conjecture on STSs.

\begin{proposition}\label{reason1}
For any positive integers $j$ and $d$, any $(j+3,j+d)$-configuration in a $4${\upshape CS} contains
a $(j+3,j)$-configuration as a substructure.
\end{proposition}

\noindent
{\bf Proof.}\quad
If a $(j+3,j+d)$-configuration contains a $4$-cycle, say $C$, in which
each vertex is also contained in another $4$-cycle,
then by discarding $C$ we
obtain a $(j+3,j+d-1)$-configuration. We prove that
for any positive integer $d$ a $(j+3,j+d)$-configuration contains
such a $4$-cycle. Suppose to the contrary that each $4$-cycle in
a given $(j+3,j+d)$-configuration ${\mathcal A}$ has at least one vertex appearing in
no other $4$-cycles.
If $d \geq 4$, the total number of vertices exceeds $j+3$, a contradiction.
Hence, we have $d=1$, $2$ or $3$.
However, by counting the total number of vertices,
it is easy to see that each case yields a contradiction.
\hfill $\Box$

\bigskip

Proposition \ref{reason1} says that any denser configuration on
$j+3$ vertices, including a $4$CS or a maximum packing, contains a
$(j+3,j)$-configuration as its substructure. On the other hand, for
$j=2$ and $1\leq e \leq 3$ every nontrivial $4$CS$(v)$ contains
$(j+3+e,j)$-configurations (see Bryant et al.\ \cite{BGGM}).
However, as we will see in the next section, we can construct a
$4$CS containing no $(j+3,j)$-configurations for any $j$ satisfying
$2 \leq j \leq 4$. Therefore, it may be natural to ask the following
question similar to Erd\H{o}s' conjecture:

\begin{problem}\label{sparse4CS}
Does there exist for every integer $r \geq 3$ a constant number $v_0(r)$
such that if $v > v_0(r)$ and $v$ is admissible, then there exists a $4${\upshape CS}$(v)$
containing no $(j+3,j)$-configurations for any $j$ satisfying $2\leq j\leq r$?
\end{problem}

\noindent
{\bf Remark.}\quad
While for any positive integers $e$ and $j$
every nontrivial STS on a sufficiently large number of vertices
contains a $(j+2+e,j)$-configuration,
we do not know in general the behavior of $(j+3+e,j)$-configurations except for $j=2$.
We briefly discuss in Section $3$ the maximum number of $4$-cycles
of a $4$-cycle packing avoiding $(j+3,j)$-configurations.

Following the terminology of STSs, we say that a $4$CS is {\it $r$-sparse} if it contains no
$(j+3,j)$-configuration for any $j$ satisfying $2\leq j\leq r$.
Every $r$-sparse $4$CS is also $(r-1)$-sparse for $r \geq 3$.
Since no $(5,2)$-configuration can appear in a $4$CS, every $4$CS is $2$-sparse.
Up to isomorphism, there are two kinds of $(6,3)$-configuration described by three $4$-cycles
$(a,b,c,d)$, $(a,e,c,f)$ and $(b,e,f,d)$, and $(a,b,c,d)$, $(a,e,c,f)$
and $(b,e,d,f)$ respectively.
A routine argument proves that any $(7,4)$-configuration is isomorphic
and can be described by four $4$-cycles
$(a,b,c,d)$, $(a,e,b,f)$, $(c,f,d,g)$ and $(a,c,e,g)$.
Hence, a $4$CS is $3$-sparse if it lacks the two types of $(6,3)$-configuration, and
it is $4$-sparse if it also avoids the unique type of $(7,4)$-configuration simultaneously.

Our results presented in the next section give resolution for the existence
problem of a $4$-sparse $4$CS$(v)$.

\begin{theorem}\label{4sparse}
There exists a $4$-sparse $4${\upshape CS}$(v)$
if and only if $v \equiv 1$ {\rm (mod $8$)}.
\end{theorem}

Up to isomorphism, there are four possible configurations formed by two $4$-cycles
in a $4$CS, the numbers of vertices ranging from six to eight.
While there are two kinds of $(6,2)$-configuration, both $(7,2)$- and $(8,2)$-configurations
are unique.
A $(6,2)$-configuration sharing a common diagonal,
described by two $4$-cycles $(a,b,c,d)$ and $(a,e,c,f)$,
is called the {\it double-diamond} configuration.
A $4$-cycle system is said to be {\it $D$-avoiding}
if it contains no double-diamond configurations.

Bryant et al.\ \cite{BGGM} showed that
for every admissible order $v$ there exists a $D$-avoiding $4$CS$(v)$.

\begin{theorem}\label{D-avoiding}{\rm\bf{(Bryant et al.) \cite{BGGM}}}
There exists a $D$-avoiding $4${\upshape CS}$(v)$ for all $v \equiv 1$ {\rm (mod $8$)}.
\end{theorem}

Since a double-diamond configuration appears in both types of $(6,3)$-configuration,
every $D$-avoiding $4$CS is $3$-sparse but the converse does not hold.
In fact, for every small admissible order $v$ one can easily find a $3$-sparse $4$CS$(v)$ which is not $D$-avoiding.
On the other hand, Bryant et al.\ \cite{BGGM} showed that
the other type of $(6,2)$-configuration appears
constantly depending only on the order $v$, that is,
the number of occurrences is unique between $4$CSs
of the same order.
Considering these facts,
we say that a $4$CS is {\it strictly} $r$-sparse if it is both $r$-sparse and
$D$-avoiding.

In Section $2$, we give a proof of existence of a strictly $4$-sparse $4$CS$(v)$
for every admissible order $v$.

\begin{theorem}\label{strictly4}
There exists a strictly $4$-sparse $4${\upshape CS}$(v)$
if and only if $v \equiv 1$ {\rm (mod $8$)}.
\end{theorem}

We also study in Section $3$ the maximum number of $4$-cycles of a $4$-cycle packing
avoiding $(j+3,j)$-configurations.

Let $ex(v,r)$ be the maximum number of $4$-cycles of a $4$-cycle packing
of order $v$ containing neither double-diamond configurations
nor $(j+3,j)$-configurations for every $2 \leq j \leq r$.
By probabilistic methods, we prove that for any positive integer $r \geq 2$
the maximum number $ex(v,r) = O(v^2)$.

\section{Strictly $4$-sparse $4$-cycle systems}
In this section, we present a proof of Theorem \ref{strictly4}.
Obviously, the proof also verifies Theorem \ref{4sparse}. To show
Theorem \ref{strictly4}, we first prove two lemmas.

A {\it jointed-diamond} configuration in a $4$CS is a $(7,3)$-configuration
described by three $4$-cycles $(a,b,c,d)$, $(a,e,b,g)$ and $(c,f,d,g)$;
the $4$-cycle $(a,b,c,d)$ is referred to as a {\it joint $4$-cycle}.
Every $(7,4)$-configuration contains a jointed-diamond configuration as its substructure.

\begin{lemma}\label{primepower}
Let $q$ be a prime power satisfying $q \equiv 1$ {\rm (mod $8$)} and
not a power of three.
Then there exists a strictly $4$-sparse $4${\upshape CS}$(q)$.
\end{lemma}

\noindent {\bf Proof.}\quad Let $q$ be a prime power satisfying $q
\equiv 1$ (mod $8$) and not a power of three. Let $\chi$ be a
multiplicative character of order four of GF$(q)$ such that
$\chi(x)$ has possible values $1$, $-1$, $i$, $-i$ for $x \not=0$.
Then there exists a $4$-cycle $(0,x,x-1,x^2)$, $x \in {\rm GF}(q)$,
such that $\chi(x^2)=-1$, $\chi((x^2-x+1)^2)=-1$, and
$\chi(x(x^2-x+1))) = 1$ (see Bryant et al.\ \cite{BGGM}).
Considering these conditions, we have either $\chi(x)=i$,
$\chi(x^2-x+1)=-i$, and $\chi(x(x-1))=i\cdot\chi(x-1)$, or
$\chi(x)=-i$, $\chi(x^2-x+1)=i$, and
$\chi(x(x-1))=-i\cdot\chi(x-1)$. Also, since $q \equiv 1$ (mod $8$),
we have $\chi(-1)=1$. Let $\alpha$ be a primitive element of GF$(q)$
and $V$ the set of all elements of GF$(q)$. Define a set ${\mathcal
C}$ of $4$-cycles as
$\{y,x\cdot\alpha^{4n}+y,(x-1)\cdot\alpha^{4n}+y,x^2\cdot\alpha^{4n}+y:
y \in {\rm GF}(q), 0\leq n \leq \frac{q-1}{8}-1\}$. Then
$(V,{\mathcal C})$ forms a $D$-avoiding $4$CS$(q)$. In fact,
${\mathcal C}$ is developed from the $4$-cycle $(0,x,x-1,x^2)$ by
the group $G=\{z \mapsto z\cdot \alpha^{4n}+y: y,z\in {\rm GF}(q),
0\leq n \leq \frac{q-1}{8}-1\}$. To prove that $(V,{\mathcal C})$ is
strictly $4$-sparse, it suffices to show that $(V,{\mathcal C})$
contains no jointed-diamond configurations. Suppose to the contrary
that it contains a jointed-diamond configuration $J$ described by
three $4$-cycles $(a,b,c,d)$, $(a,e,b,g)$ and $(c,f,d,g)$. Since
every $4$-cycle in ${\mathcal C}$ can be obtained from
$(0,x,x-1,x^2)$ by the group $G$, considering the joint $4$-cycle
$(a,b,c,d)$, we have $\chi(a-b)=-\chi(c-d)$. However, since the
edges $\{a,b\}$ and $\{c,d\}$ lie in diagonals of $(a,e,b,g)$ and
$(c,f,d,g)$ respectively, we have $\chi(a-b)=\chi(c-d)$,
$i\cdot\chi(c-d)$ or $-i\cdot\chi(c-d)$, a contradiction. The proof
is complete. \hfill $\Box$

\begin{lemma}\label{order9}
There exists a strictly $4$-sparse $4${\upshape CS}$(9)$.
\end{lemma}

\noindent
{\bf Proof.}\quad
Let $V = \{0,1,2,\dots,8\}$ be the set of elements of the cyclic group
${\textit{\textbf{Z}}}_{9}$. Define a set ${\mathcal C}$ of $4$-cycles as
$\{(0+a,1+a,8+a,5+a): a \in {\textit{\textbf{Z}}}_{9}\}$.
The pair $(V,{\mathcal C})$ forms a $4$CS$(9)$ under the transitive action
of ${\textit{\textbf{Z}}}_{9}$ on the vertex set $V$. Since ${\mathcal C}$ has only one
$4$-cycle orbit, $(V,{\mathcal C})$ is $D$-avoiding, and hence it is $3$-sparse.

Suppose to the contrary that $(V,{\mathcal C})$ is not $4$-sparse and
contains a jointed-diamond.
Take a representative, say $C=(0,1,8,5)$, of the $4$-cycle orbit.
The two differences of the vertices in a diagonal of $C$
are $\pm 1$ and $\mp 4$ respectively.
Hence, the joint $4$-cycle in a jointed-diamond lying in ${\mathcal C}$ has the form
$(a,b,c,d)$, where the differences $a-b$ and $c-d$ are each $1$, $-1$, $4$ or $-4$.
However, considering the four differences of the adjacent vertices in $C$,
this is a contradiction.\\
\hfill $\Box$

\bigskip

We now return to the proof of Theorem \ref{strictly4}.
The proof employs a special decomposition of
the complete graph into smaller complete graphs.

A {\it group divisible design} with {\it index} one is a triple
$(V, {\mathcal G}, {\mathcal B})$, where
\begin{enumerate}
    \item[(i)] $V$ is a finite set of elements called points,
        \item[(ii)] ${\mathcal G}$ is a family of subsets of $V$, called {\it groups},
                which partition $V$,
        \item[(iii)] ${\mathcal B}$ is a collection of subsets of $V$, called
                {\it blocks}, such that every pair of points from distinct
                groups occurs in exactly one blocks,
        \item[(iv)] $|G \cap B| \leq 1$ for all $G \in {\mathcal G}$ and $B
                \in {\mathcal B}$.
\end{enumerate}
When all blocks are of the same size $k$ and the number of groups of
size $n_i$ is $t_i$, one refers to the design as a $k$-GDD of {\it
type} $n_0^{t_0}n_1^{t_1}\cdots n_{g-1}^{t_{g-1}}$, where
$t_0+t_1+\dots+t_{g-1}= |{\mathcal G}|$. We need $4$-GDDs and the
required types are of $12^t$ $(t \geq 4)$, $4^{3t+1}$ $(t \geq 1)$,
$8^{3t+1}$ $(t \geq 1)$, and $2^{3t}5^1$ $(t \geq 3)$. For their
existence, we refer the reader to Colbourn and Dinitz \cite{CD}. \\
\\
{\bf Proof of Theorem \ref{strictly4}.}\quad A strictly $4$-sparse
$4$CS$(v)$ is necessarily $D$-avoiding. We follow a part of the
proof of existence of a $D$-avoiding $4$CS$(v)$ by Bryant et al.\
\cite{BGGM} and consider four cases:

\medskip
{\it Case} (1) : $v \equiv 1$ (mod $24$).
Lemma \ref{primepower} gives a strictly $4$-sparse $4$CS$(v)$ for
$v \leq 73$ and $v \equiv 1$ (mod $24$). We consider the case $v > 73$.
Take a $4$-GDD $(V,{\mathcal B}, {\mathcal G})$ of type $12^t$ for $t \geq 4$.
For each group $G \in {\mathcal G}$, take $(G\times \{0,1\}) \cup \{\infty\}$
by replacing each point by two new points and adding a new point $\infty$.
Let ${\mathcal H}_G$ be a copy of the strictly $4$-sparse $4$CS$(25)$
given in Lemma \ref{primepower} on $(G\times \{0,1\}) \cup \{\infty\}$.
For each block $B=\{a,b,c,d\} \in {\mathcal B}$, construct a $4$-cycle decomposition
${\mathcal C}_B$ of a copy of $K_{2,2,2,2}$ on $B \times \{0,1\}$
by developing a $4$-cycle $((a,0),(b,0),(c,1),(d,0))$ under the group
$\langle (d)(a\ b\ c)\rangle\times {\textit{\textbf{Z}}}_2$.
Let $W = (V\times \{0,1\})\cup\{\infty\}$ and
${\mathcal D} = (\bigcup_{G\in{\mathcal G}} {\mathcal H}_G)
\cup (\bigcup_{B\in{\mathcal B}}{\mathcal C}_B)$.
Then $(W,{\mathcal D})$ forms a $4$CS$(24t+1)$.
Since no pair of $4$-cycles in ${\mathcal D}$ shares a common
diagonal, $(W,{\mathcal D})$ is $D$-avoiding.

It remains to establish that the $4$CS contains no $(7,4)$-configuration.
Suppose to the contrary that  $(W,{\mathcal D})$ contains a  $(7,4)$-configuration.
Then it contains a jointed-diamond configuration $J$.
If the joint $4$-cycle in $J$ lies in ${\mathcal H}_G$,
the other two $4$-cycles in $J$
are also in ${\mathcal H}_G$. Since ${\mathcal H}_G$ is a copy of
a strictly $4$-sparse $4$CS$(25)$, this is a contradiction.
If the joint $4$-cycle in $J$ lies in ${\mathcal C}_B$, again
the other two $4$-cycles in $J$ are in ${\mathcal C}_B$.
A routine argument proves that ${\mathcal C}_B$ contains no
jointed-diamond configuration.

\medskip
{\it Case} (2) : $v \equiv 9$ (mod $24$).
Lemma \ref{order9} gives a strictly $4$-sparse $4$CS$(9)$.
Take a $4$-GDD $(V,{\mathcal B}, {\mathcal G})$ of type $4^{3t+1}$ for $t \geq 1$.
As in {\it Case} (1), construct a $4$CS$(24t+9)$ on $(V \times \{0,1\})\cup\{\infty\}$
by placing a copy of the strictly $4$-sparse $4$CS$(9)$ given in Lemma \ref{order9} and
decomposing $K_{2,2,2,2}$s into $4$-cycles. By following the argument in {\it Case} (1),
the resulting $4$CS$(24t+9)$ is strictly $4$-sparse.

\medskip
{\it Case} (3) : $v \equiv 17$ (mod $48$).
Employing the strictly $4$-sparse $4$CS$(17)$ constructed in Lemma \ref{primepower}
and a $4$-GDD of type $8^{3t+1}$ for $t \geq 1$, we obtain the required
strictly $4$-sparse $4$CSs by the same technique as in {\it Case} (1).

\medskip
{\it Case} (4) : $v \equiv 41$ (mod $48$). Lemma \ref{primepower}
gives a strictly $4$-sparse $4$CS$(v)$ for $v \leq 137$ and $v
\equiv 41$ (mod $48$). We consider the case $v > 137$. Take a
$4$-GDD $(V,{\mathcal B}, {\mathcal G})$ of type $2^{3t}5^1$ for $t
\geq 3$. For each block $B=\{a,b,c,d\} \in {\mathcal B}$, replace
each point in $B$ by four new points and define $A_i = \{i\}\times
\{0,1,2,3\}$ for $i \in B$. The points and lines of an affine space
over GF$(2^2)$ of dimension $2$ form a $4$-GDD of type $4^4$. For
each $B \in {\mathcal B}$, place a $4$-GDD of type $4^4$ on $B\times
\{0,1,2,3\}$ such that the set of groups is $\{A_i: i \in B\}$ and
let ${\mathcal C}_B$ be the resulting blocks of the $4$-GDD on
$B\times \{0,1,2,3\}$. For each ${\mathcal C}_B$, $B \in {\mathcal
B}$, construct a $4$-cycle decomposition ${\mathcal D}_{{\mathcal
C}_B}$ of a copy of $K_{2,2,2,2}$ on ${\mathcal C}_B \times \{0,1\}$
by developing a $4$-cycle $((a,i,0),(b,j,0),(c,k,1),(d,l,0))$ under
the group $\langle ((d,l))((a,i)\ (b,j)\ (c,k))\rangle\times
{\textit{\textbf{Z}}}_2$. For each group $G \in {\mathcal G}$, take
$(G\times \{0,1,2,3\}\times\{0,1\}) \cup \{\infty\}$ and let
${\mathcal H}_G$ be a copy of either the strictly $4$-sparse
$4$CS$(17)$ or $4$CS$(41)$ given in Lemma \ref{primepower} on
$(G\times \{0,1,\dots,7\}) \cup \{\infty\}$ according to the group
size $|G|$, that is, place a copy of the $4$CS$(17)$ if $|G|=2$,
otherwise put a copy of the $4$CS$(41)$. Let $W = (V\times
\{0,1,2,3\}\times \{0,1\})\cup\{\infty\}$ and ${\mathcal
E}=(\bigcup_{G\in{\mathcal G}} {\mathcal H}_G) \cup
(\bigcup_{B\in{\mathcal B}}{\mathcal D}_{{\mathcal C}_B})$. It is
straightforward to see that $(W,{\mathcal E})$ forms a
$4$CS$(48t+41)$. The same argument as in {\it Case} (1) proves that
$(W,{\mathcal E})$ is strictly $4$-sparse. \hfill $\Box$

\section{$r$-Sparse $4$-cycle packing}
In this section, we consider the maximum number of $4$-cycles in a
$4$-cycle packing of order $v$ avoiding $(j+3,j)$-configurations.
As with a $4$CS, a $4$-cycle packing is said to be {\it $r$-sparse}
if it contains no $(j+3,j)$-configuration for any $j$
satisfying $2\leq j\leq r$. Also if it is $r$-sparse and $D$-avoiding,
we say that it is {\it strictly} $r$-sparse.
We prove that for any positive integer $r \geq 2$ and sufficiently large
integer $v$ there exists a constant number $c$ such that there exists a
strictly $r$-sparse $4$-cycle packing of order $v$ with $c\cdot v^2$ $4$-cycles.
It is notable that a resolution for the analogous problem to the $r$-sparse
conjecture on STSs would prove that $c \sim \frac{1}{8}$.

Let ${\mathcal F}$ be a set of configurations
of $4$-cycles and $ex(v,{\mathcal F})$ the largest positive integer $n$ such that
there exists a set ${\mathcal C}$ of $n$ $4$-cycles on a finite set $V$ of cardinality
$v$ having property that ${\mathcal C}$ contains no configuration which is isomorphic
to a member $F \in {\mathcal F}$.

\begin{theorem}\label{bound}
For any positive integer $r \geq 2$ and sufficiently large
integer $v$ there exists a constant number $c$ such that
there exists a strictly $r$-sparse $4$-cycle packing of order $v$ with
$c\cdot v^2$ $4$-cycles.
\end{theorem}

\noindent
{\bf Proof.}\quad
Let $V$ be a finite set of cardinality $v$. Define ${\mathcal F}'$ as the set of
all nonisomorphic $(j+3,j)$-configurations for $2\leq j \leq r$ and ${\mathcal F}''$
as the set of all nonisomorphic $(4,2)$- and $(6,2)$-configurations.
Let ${\mathcal F} = {\mathcal F}' \cup {\mathcal F}''$. It is easy to see that
if $ex(v,{\mathcal F}) \geq c\cdot v^2$ for some constant $c$,
then the assertion of Theorem \ref{bound} follows.

Pick uniformly at random $4$-cycles from $V$ with probability $p=\frac{c'}{v^2}$,
independently of the others, where $c'$ satisfies $0<c'<\frac{1}{44}$.
Let $b_C$ be a random variable counting the number of configurations isomorphic to $C$ in the resulting
set of $4$-cycles and $E(b_C)$ its expected value. Then

\begin{eqnarray*}
E\left( \sum_{C\simeq F\in{\mathcal F}}b_C \right) &\leq&
{{v}\choose{4}}\cdot{3\cdot{{4}\choose{4}}\choose{2}}\cdot p^2+
{{v}\choose{6}}\cdot{3\cdot{{6}\choose{4}}\choose{2}}\cdot p^2\\
&&{+}\: \sum_{j=2}^{r}{{v}\choose{j+3}}\cdot{3\cdot{{j+3}\choose{4}}\choose{j}}\cdot p^{j}\\
&\leq&\left[ {{v}\choose{4}}\cdot{{3}\choose{2}}+{{v}\choose{6}}\cdot{{45}\choose{2}}\right]p^2\\
&&{+}\: \sum_{j=2}^{r} \left( \frac{e\cdot v}{j+3} \right)^{j+3}
\cdot \left( \frac{e\cdot (j+3)^3}{8} \right)^j \cdot p^j\\
&=&\frac{11\cdot c'^2}{8}\cdot v^2+f(v),
\end{eqnarray*}
where $f(v)=O(v)$. By Markov's Inequality,
\[P\left( \sum_{C\simeq F\in{\mathcal F}}b_C \geq 2\cdot E\left( \sum_{C\simeq F\in{\mathcal F}}b_C \right)\right)
\leq \frac{1}{2}.\]
Hence,
\[P\left( \sum_{C\simeq F\in{\mathcal F}}b_C \leq \frac{11\cdot c'^2}{4}\cdot v^2+2\cdot f(v)\right)
\geq \frac{1}{2}.\]

Let $t$ be a random variable counting the number of $4$-cycles and $E(t)$ its expected value.
Then
\[E(t)=p\cdot3\cdot{{v}\choose{4}}=\frac{c'}{8}\cdot v^2 - g(v),\]
where $g(v)=O(v)$. Since $t$ is a binomial random variable,
we have for sufficiently large $v$
\begin{eqnarray*}
P\left( t<\frac{E(t)}{2}\right) &<& e^{-\frac{E(t)}{8}}<\frac{1}{2}.
\end{eqnarray*}

Hence, if $v$ is sufficiently large, then we have, with positive probability,
a set ${\mathcal S}$ of $4$-cycles with the property
that $|{\mathcal S}| \geq \frac{E(t)}{2}$ and
the number of configurations in ${\mathcal S}$ which are isomorphic
to a member of ${\mathcal F}$ is at most $\frac{11\cdot c'^2}{4}\cdot v^2+2\cdot f(v)$.
Since $f(v),\ g(v)=O(v)$, by deleting a $4$-cycle from each configuration
isomorphic to a member of ${\mathcal F}$, we have
\[ex(v,{\mathcal F}) \geq \frac{c'(1-44\cdot c')}{16}\cdot v^2 - h(v),\]
where $h(v)=O(v)$.
The proof is complete.
\hfill $\Box$

\section*{Acknowledgements}
The authors would like to thank Professor Masakazu Jimbo for helpful
comments. Research of the first author is supported by JSPS Research
Fellowships for Young Scientists.

\end{document}